# μ- Integrable Functions and Weak Convergence of Finite Measures

**Renying Zeng**

Mathematics Department, Saskatchewan Polytechnic, 1130 Idylwyld Dr. N, Saskatoon, SK Canada S7L 4J7

Email: renying.zeng@saskpolytech.ca

**Abstract** This paper deals with functions that defined in metric spaces and valued in complete paranormed vector spaces or valued in Banach spaces, and obtains some necessary and sufficient conditions for weak convergence of finite measures. Our main results are as follows. Let $X$ be a Banach space, $(\Omega, \rho)$ be a compact metric space, $\Sigma$ the σ-algebra generated by all open sets—the smallest σ-ring containing all the open sets. A sequence $\{\mu_n\}$ of probability measure is weakly convergent to probability measure $\mu$ if and only if $\lim_{n\to\infty}\int_\Omega g(s)d\mu_n = \int_\Omega g(s)d\mu$ for every bounded continuous function $g: \Omega \to X$. If $X$ is a Frechet space in Bourbaki's terminology, and has a base, then $\{\mu_n\}$ is weakly convergent to $\mu$ if and only if $\lim_{n\to\infty}\int_\Omega g(s)d\mu_n = \int_\Omega g(s)d\mu$ for every bounded continuous function $g: \Omega \to X$.

## 1. Introduction and Terminology

If $(\Omega, \Sigma)$ is a measurable space, a sequence $\{\mu_n\}$ of probability measure is called weakly convergent to probability measure $\mu$ if $\lim_{n\to\infty}\int_\Omega g(s)d\mu_n = \int_\Omega g(s)d\mu$ for every bounded continuous function $g: \Omega \to (-\infty, +\infty)$. Nielsen [1] proved that, if $X$ is a smooth Banach space, we have $\lim_{n\to\infty}\int_\Omega g(s)d\mu_n = \int_\Omega g(s)d\mu$ for every bounded continuous function $g: \Omega \to X$. Yang [2] discussed the situation of $g: \Omega \to (-\infty, +\infty)$, and proved that the function $g$ is not necessarily to be point-wise continuous. Wei [3] worked with functions valued in a metric space $X$, and assumed that $\{\mu_n\}$ is tight and weak convergence of the finite dimensional distributions of $\{\mu_n\}$ to $\mu$.

In this paper, we extended the necessary and sufficient condition to any Banach space, and to any Frechet space (in Bourbaki's terminology) that has a base. We generalized Nielsen's result in [1], and generalized the corresponding results of Yang and Wei in [2, 3].

Let $K$ be the field of real numbers or the field of complex numbers, and $X$ be a vector space over $K$. A **paranormed space** is a pair $(X, \|\cdot\|)$, where $\|\cdot\|$ is a function, called a **paranorm**, such that

(a) $\|x\| \geq 0$, $\|x\| = 0 \Leftrightarrow x = 0$;
(b) $\|x+y\| \leq \|x\| + \|y\|$;
(c) $\|-x\| = \|x\|$;
(d) $\lim_{\alpha \to 0} \|\alpha x\| = 0$, $\lim_{x \to 0} \|\alpha x\| = 0$.

Since $\|x - y\| \leq \|x\| + \|y\|$, for $\forall x, y \in X$, $\|x - y\|$ defines a **metric** in a paranormed space.

In what follows paranormed paces will always be regarded as metric spaces with respect to the metric $\|\cdot\|$.

It is known that a normed vector space is a paranormed vector space, but a paranormed space is not necessary to be a normed vector space.

A sequence $\{x_n\}$ in $X$ is said to be a **Cauchy sequence** if $\forall \varepsilon > 0$ $\exists N > 0$, when $m, n > N$, $\|x_m - x_n\| < \varepsilon$. A paranormed vector space is said to be **complete** if every Cauchy sequence $\{x_n\}$ is convergent in $X$.

Any Banach space is a complete paranormed space. But the converse is not true.

A complete paranormed space is called a **Frechet space** in Bourbaki's terminology, see [4].

Let $\Omega$ be an arbitrary non-empty set. A family $\Sigma$ of subsets of $\Omega$ is called **ring** if

(a) $\Sigma$ contains **the empty set** $\Phi$;
(b) if $A, B \in \Sigma$ then $A \cup B \in \Sigma$ and $A \setminus B \in \Sigma$.

It follows that also the intersection $A \cap B$ belongs to $\Sigma$ because

$$A \cap B = B \setminus (B \setminus A)$$

is also in $\Sigma$. Hence, a ring is closed under taking the set-theoretic operations $\cap, \cup, \setminus$.

A ring $\Sigma$ is called a **σ-ring** if the union $\bigcup_{k=1}^{\infty} A_k$ of any countable family $\{A_k\}_{k=1}^{\infty}$ of sets from $\Sigma$ is also in $\Sigma$. A σ-ring in $\Omega$ that contains $\Omega$ is called a **σ-algebra** of $\Omega$.

Let $\Omega$ be an arbitrary non-empty set and $\Sigma$ be a σ-algebra of $\Omega$, R+ the set of all non-negative real numbers. A functional μ: $\Sigma \to$ R+ called a **σ-additive measure** on $\Sigma$ if whenever a set $A \in \Sigma$ is a disjoint union of an at most countable sequence in $\Sigma$, i.e., $A = \bigcup_{k=1}^{N} A_k$ (where $N$ is either finite or $N = \infty$) and $A_k \cap A_l = \Phi, k \neq i$, then

$$\mu(A) = \sum_{k=1}^{N} \mu A_k.$$

If $N = \infty$ then the above sum is understood as a series.

A subset $A \subset \Omega$ is called **μ-measurable** if $A \in \Sigma$.

A **σ-additive measure** $\mu$ on $\Sigma$ is called a **finite measure** if $\Omega \in \Sigma$ and $\mu(\Omega) < \infty$. $(\Omega, \Sigma, \mu)$ is called a **finite measure space** if the measure $\mu$ on $\Sigma$ is finite.

A function $g: \Omega \to X$ is called **μ-measurable function** if for any scalar $\alpha$

$$\{s \in \Omega; g(s) < \alpha\}$$

is a **μ-measurable** subset.

A function $g: \Omega \to X$ is called a **simple function** if there exist measurable sets $B_j \in \Sigma$ with $B_i \cap B_j = \Phi \, (i \neq j)$ and $x_j \in X$ $(j = 1, 2, \ldots, k)$ such that

$$g(s) = \begin{cases} x_j, & \text{if} \quad s \in B_j \, (j = 1, \cdots, k) \\ 0, & \text{otherwise} \end{cases}$$

The **μ- Integral** of the simple function $g$ is defined as

$$\int_{\Omega} g(s) d\mu = \sum_{j=1}^{k} x_j \mu B_j.$$

A function $g: \Omega \to X$ is called **μ- Integrable**, if there exists a sequence of $\{g_n\}$ of simple functions such that

(a) $\lim_{n \to \infty} g_n(s) = g(s)$, **μ-a.e.** on $\Omega$ \hfill (1)

i.e., there exists $A \in \Sigma$ and $\mu(A) = 0$ such that $\lim_{n \to \infty} g_n(s) = g(s)$ for $s \in \Omega \setminus A$;

(b) For every continuous seminorm $p$ on $X$,

$$\lim_{n \to \infty} \int_{\Omega} p(g_n(s) - g(s)) d\mu = 0. \hfill (2)$$

In this case, $\lim_{n \to \infty} \int_{\Omega} g_n(s) d\mu$ exists, and the **μ- Integral** of $g$ is defined as

$$\int_\Omega g(s)d\mu = \lim_{n\to\infty}\int_\Omega g_n(s)d\mu.$$

When $X$ is a Banach space, the $\mu$- Integral is known as **Bochner integral** [5, 6].

A sequence $\{x_n\}$ in $X$ is called a **base** of $X$ if for each $x \in X$, there exists a unique sequence of scalars $\{\alpha_n\}$ such that

$$x = \sum_{n=1}^{\infty} \alpha_n x_n.$$

If the mappings

$$x'_n : X \to K, \quad x'_n(x) = \alpha_n$$

are continuous linear functionals on $X$, then $\{x_n\}$ is called a **Schauder base** of $X$, and one has

$$x'_i(x_j) = \begin{cases} 1, & \text{if } i = j \\ 0 & \text{if } i \neq j \end{cases}.$$

Where $x'_n$ are called the **coordinate functionals** ($n = 1, 2, \ldots$).

## 2. Some Lemmas

In what follows, let $(\Omega, \rho)$ be a compact metric space, $\Sigma$ the **σ-algebra** generated by all open sets—the smallest σ-ring containing all the open sets, (then $(\Omega, \rho, \Sigma)$ is a measurable space,) and $\Pi$ a family of σ-additive finite measures on $\Sigma$; $(X, \|\cdot\|)$ be a complete paranormed space.

Lemma 1 can be found in [7].

**Lemma 1** If a complete paranormed space $X$ has a base $\{x_n\}$, then $\{x_n\}$ is a Schauder base.

A **seminorm** is map

$$p: X \to K$$

satisfying

(a) $p(x) \geq 0$;
(b) $p(x+y) \leq p(x) + p(y)$;
(c) $p(\alpha x) = |\alpha| p(x)$ for any scalar $\alpha$.

From [4] one has the following Lemma 2.

**Lemma 2** $X$ is a complete paranormed space if and only if there is a family of continuous siminorms $P = \{p_n; n = 1, 2, ....\}$ on $X$, such that

$$p_1(x) \leq p_2(x) \leq \cdots \leq p_n(x) \leq \cdots, \quad \forall x \in X.$$

And the paranorm on $X$ can be given by

$$\|x\| = \sum_{n=1}^{\infty} \frac{1}{2^n} \frac{p_n(x)}{1 + p_n(x)}, \quad \forall x \in X.$$

Furthermore, for any topological net $\{x_\tau\} \subset X$, and $x \in X$, the following are equivalent

(a) $\lim_\tau x_\tau = x$;

(b) $\lim_\tau \| x_\tau - x \| = 0$;

(c) $\lim_\tau p(x_\tau - x) = 0$, $\forall p \in P$.

A set $A \subset X$ is called separable if $A$ has a countable dense subset, i.e., there exists a countable subset $B \subset A$ such that $\overline{B} = A$, where $\overline{B}$ is the topological closure of $B$. $B$ is called a countable dense subset of $A$.

For $A \subset X$, **span A** is the set of all possible linear combinations of the elements in $A$.

Lemma 3 has an independent interest.

**Lemma 3** Suppose $\mu \in \Pi$. A µ-measurable function $g: \Omega \to X$ is µ-integrable if and only if

(a) $g$ is µ-essential separable valued, i.e., there exists $E \in \Sigma$ with $\mu E = 0$, such that $g(\Omega \setminus E)$ is a separable subset of $X$;

(b) $\int_\Omega p(g(s))d\mu < \infty$, $\forall p \in P$, where $P$ is defined as in Lemma 2.

The following Lemma 4 is from [8].

**Lemma 4** Suppose $\{x_n\} \subset X$ such that $\overline{span\{x_n\}} = X$, with $x_n \neq 0, \forall n$. Then $\{x_n\}$ is a Schauder base in $X$ if and only if $\forall p \in P$, there exist a corresponding $q \in P$ and a constant $M > 0$, such that

$$p(\sum_{i=1}^{m} \alpha_i x_i) \leq Mq(\sum_{i=1}^{n} \alpha_i x_i)$$

for each pair of integers $m, n$ with $m \leq n$, and arbitary scalars $\alpha_1, \alpha_2, \cdots, \alpha_n$. Specially one has

$$p(\sum_{i=1}^{m} \alpha_i x_i) \leq Mq(\sum_{i=1}^{\infty} \alpha_i x_i).$$

**Lemma 5** Suppose $g_n: \Omega \to X$ ($n = 1, 2, \ldots$). If for given $p \in P$, $s_0 \in \Omega$, $\forall \varepsilon > 0$, there exists $\delta > 0$, such that $\rho(s - s_0) < \delta$ implies

$$p(g_n(s) - g_n(s_0)) < \varepsilon, \text{ for every } n,$$

then there exists $n_0$, such that $n \geq n_0$ implies

$$p(g_n(s) - g(s)) < \varepsilon \text{ for all } s \in \Omega.$$

**Lemma 6** Let $\mu, \mu_n \in \Pi$ ($n = 1, 2, \ldots$). If $\Pi$ is a relatively compact according to its metric topology, then $\forall \varepsilon > 0$, there exists a compact subset $A \subset \Omega$, such that $\mu(\Omega \setminus A) \leq \varepsilon, \mu_n(\Omega \setminus A) \leq \varepsilon$, for all $n$.

Lemma 6 can be found in [10].

**Lemma 7** If $g: \Omega \to X$ is a continuous function, then $g$ is $\mu$- integrable.

## 3. Weak Convergence of Finite Measures

The following is the traditional definition of weak convergence of measures [10].

A sequence $\{\mu_n\} \subset \Pi$ is called **weakly convergent** to $\mu \in \Pi$ if

$$\lim_{n \to \infty} \int_{\Omega} g(s) d\mu_n = \int_{\Omega} g(s) d\mu$$

for every bounded continuous function $g: \Omega \to (-\infty, +\infty)$.

This paper extends the concept of weakly convergent measures to complete paranormed spaces.

**Theorem 1** Suppose that $X$ over the number field $K$ is a complete paranormed space and has a base $\{x_n\}$ and coordinate functionals $\{x'_n\}$. Let $\mu, \mu_n \in \Pi$ ($n = 1, 2, \ldots$), then $\{\mu_n\}$ is weakly convergent to $\mu$ if and only if

$$\lim_{n\to\infty} \int_\Omega g(s)d\mu_n = \int_\Omega g(s)d\mu$$

for every bounded continuous function $g: \Omega \to X$.

Theorem 1 can be stated as the following Theorem 2, a result in probability distribution theory.

**Theorem 2** Suppose that $X$ over the number field $K$ is a complete paranormed space and has a base. Let $\varsigma_n \in \Omega$ be random elements ($n = 1, 2, \ldots$), $E$ the mathematical expectation operator, then $\{\varsigma_n\}$ converges in probability distribution to a random element $\varsigma \in \Omega$ if and only if for each bounded continuous function $g: \Omega \to X$,

$$\lim_{n\to\infty} Eg(\varsigma_n) = Eg(\varsigma).$$

If $X$ is a separable Banach space, then $X$ has a Schauder base. Similar to the proof of Theorem 1 one has the following Corollary 3.

**Corollary 3** Let $X$ be a separable Banach space. A sequence $\{\mu_n\} \subset \Pi$ is weakly convergent to $\mu \in \Pi$ if and only if

$$\lim_{n\to\infty} \int_\Omega g(s)d\mu_n = \int_\Omega g(s)d\mu$$

for every bounded continuous function $g: \Omega \to X$.

For a Banach space, however, the following Theorem 4 holds.

**Theorem 4** Let $X$ be a Banach space. A sequence $\{\mu_n\} \subset \Pi$ is weakly convergent to $\mu \in \Pi$ if and only if

$$\lim_{n\to\infty} \int_\Omega g(s)d\mu_n = \int_\Omega g(s)d\mu$$

for every bounded continuous function $g: \Omega \to X$.

In the proof of the sufficiency of Theorem 1, Replacing $X$ by $\overline{span\{Y\}}$ and replacing the seminorm $p \in P$ by the norm $\|\cdot\|$, we complete the proof of Theorem 4.

Similarly, the following Theorem 5 holds.

**Theorem 5** Suppose that $X$ over the number field $K$ is a complete paranormed space, $\overline{span\{Y\}}$ has a base, where

$$Y = g(\Omega)\cup(\cup_{n=1}^{\infty} g_n(\Omega)).$$

Let $\mu, \mu_n \in \Pi$ ($n = 1, 2, \ldots$), then $\{\mu_n\}$ is weakly convergent to $\mu$ if and only if

$$\lim_{n\to\infty}\int_{\Omega} g(s)d\mu_n = \int_{\Omega} g(s)d\mu$$

for every bounded continuous function $g: \Omega \to X$.

Theorem 4 and Theorem 5 can be rewritten as the following Theorem 6 and Theorem 7.

**Theorem 6** Suppose that $X$ is a Banach space. Let $\varsigma_n \in \Omega$ be random elements ($n = 1, 2, \ldots$), $E$ the mathematical expectation operator, then $\{\varsigma_n\}$ converges in probability distribution to a random element $\varsigma \in \Omega$ if and only if for each bounded continuous function $g: \Omega \to X$,

$$\lim_{n\to\infty} Eg(\varsigma_n) = Eg(\varsigma).$$

**Theorem 7** Suppose that $X$ over the number field $K$ is a complete paranormed space, $\overline{span\{Y\}}$ has a base, where $Y = g(\Omega)\cup(\cup_{n=1}^{\infty} g_n(\Omega))$. Let $\varsigma_n \in \Omega$ be random elements ($n = 1, 2, \ldots$), $E$ the mathematical expectation operator, then $\{\varsigma_n\}$ converges in probability distribution to a random element $\varsigma \in \Omega$ if and only if for each bounded continuous function $g: \Omega \to X$,

$$\lim_{n\to\infty} Eg(\varsigma_n) = Eg(\varsigma).$$

### 4. Conclusion Remark

Let $\Omega$ be a non-empty set and $\Sigma$ be a σ-algebra of $\Omega$. A measure $\mu: \Sigma \to [0,+\infty]$ is called a finite measure if

$$\mu(\Omega) < \infty;$$

a probability measure if

$$\mu(\Omega) = 1.$$

This paper obtains some necessary and sufficient conditions for the weak convergence of finite measures in complete paranormed vector spaces that have bases and in general Banach spaces. Or, actually, we extend the concept of weak convergence of probability measures to complete Paranormed vector spaces that have bases, and to general Banach spaces.

We generalize also the concept of convergence of random variables in probability distributions, to Paranormed vector spaces and to general Banach spaces.